\documentclass{amsart}

\usepackage{amsmath}
\usepackage{amsthm}
\usepackage{amsfonts}
\usepackage{amssymb}
\usepackage{latexsym}
\usepackage{amscd}
\usepackage{stmaryrd}

\theoremstyle{plain}
\newtheorem{thm}{Theorem}[section]

\newtheorem{lem}[thm]{Lemma}
\newtheorem{prop}[thm]{Proposition}

\theoremstyle{definition}
\newtheorem{rem}{Remark}[section]

\def\proof{\noindent {\bf Proof.\;}}

\numberwithin{equation}{section}

\begin{document}

\title{On functional relations for the alternating analogues of Tornheim's double zeta function}

\date{}

\author{Zhong-hua Li}

\maketitle

\begin{center}
{\small Department of Mathematics, Tongji University, No. 1239 Siping Road},\\
{\small Shanghai 200092, China}\\
{\small E-mail address: zhonghua\_li@tongji.edu.cn}
\end{center}

\vskip10pt

{\footnotesize
\begin{quote}

\noindent {\bf Abstract.}
We give new proofs of two functional relations for the alternating analogues of Tornheim's double zeta function. Using the functional relations, we give new proofs of some evaluation formulas found by H. Tsumura for these alternating series.


\noindent{\bf 2010MSC}: 11M32, 40B05
\end{quote}
}


\section{Introduction}

Let $\mathbb{N}$ be the set of positive integers, and $\mathbb{C}$ the field of complex numbers.

For $s_1,s_2,s_3\in\mathbb{C}$ with $\Re(s_1+s_3)>1$, $\Re(s_2+s_3)>1$ and $\Re(s_1+s_2+s_3)>2$,
Tornheim's double zeta function is defined as
\begin{equation}
T(s_1,s_2,s_3):=\sum\limits_{m,n=1}^\infty\frac{1}{m^{s_1}n^{s_2}(m+n)^{s_3}}.
\label{Eq:defn-Tornheim}
\end{equation}
It has two alternating analogues:
\begin{equation}
R(s_1,s_2,s_3):=\sum\limits_{m,n=1}^\infty\frac{(-1)^n}{m^{s_1}n^{s_2}(m+n)^{s_3}},
\label{Eq:defn-R}
\end{equation}
and
\begin{equation}
S(s_1,s_2,s_3):=\sum\limits_{m,n=1}^\infty\frac{(-1)^{m+n}}{m^{s_1}n^{s_2}(m+n)^{s_3}}.
\label{Eq:defn-S}
\end{equation}
It is proved in \cite[Theorem 1]{matsumoto02} that the function $T(s_1,s_2,s_3)$ can be meromorphically continued to the whole $\mathbb{C}^3$-space, and its singularities are the subsets of $\mathbb{C}^3$ defined by one of the equations
$$s_1+s_3=1-l,s_2+s_3=1-l (l\in\mathbb{N}\cup\{0\}),s_1+s_2+s_3=2.$$
While from \cite[Theorem 2.1]{nakamura081}, we know that $R(s_1,s_2,s_3)$ can be meromorphically continued to $\mathbb{C}^3$ with the singularities lying on the subsets of $\mathbb{C}^3$ defined by one of the equations $s_1+s_3=1-l$ ($l\in\mathbb{N}\cup\{0\}$), and $S(s_1,s_2,s_3)$ can be analytically continued to $\mathbb{C}^3$.

The function $T(s_1,s_2,s_3)$ is a generalization of the Riemann zeta function $\zeta(s)$. As the Riemann zeta value $\zeta(k)$ ($k\in \mathbb{N}, k>1$) is important in the study of $\zeta(s)$, it is natural to consider the values of $T(s_1,s_2,s_3)$ with all arguments are positive integers. In fact, L. Tornheim \cite{tornheim} first introduced the series $T(p,q,r)$ for $p,q,r\in\mathbb{N}$ in 1950's. And since then,
a lot of results on evaluating the values $T(p,q,r)$ in terms of Riemann zeta values have been found. See for example
\cite{huard,mordell,nakamura,subbarao,tornheim,tsumura2002,tsumura2007} and the references therein.
In \cite{tsumura2007}, by Fourier expansion technique, H. Tsumura proved a functional relation for Tornheim's double zeta function which represents
$T(a,b,s)+(-1)^bT(b,s,a)+(-1)^aT(s,a,b)$ with $a,b\in\mathbb{N}$ and $s\in\mathbb{C}$ via Riemann zeta function.
After that, by introducing the function
$$\sum\limits_{m\neq 0,n\neq 0\atop m+n\neq 0}\frac{1}{m^{s_1}n^{s_2}(m+n)^{s_3}},$$
T. Nakamura \cite{nakamura} gave a \lq\lq simpler" version of this functional relation.
And in \cite{matsumoto}, K. Matsumoto, T. Nakamura, H. Ochiai and H. Tsumura showed that these two functional relations are the same.

The alternating analogues of Tornheim's double zeta series were first introduced by M. V. Subbarao and R. Sitaramachandrarao in
\cite{subbarao}. They posed the problem to evaluate $S(r,r,r)$ and $R(r,r,r)$ for $r\in\mathbb{N}$. In a series of papers
\cite{tsumura2002,tsumura2003,tsumura2004,tsumura2009}, H. Tsumura obtained some fascinating results on evaluating $S(p,q,r)$ and
$R(p,q,r)$ by using Fourier expansion technique. He gave an evaluation formula for $S(r,r,r)$ for any positive odd integer $r$ in \cite{tsumura2002}, and for
$R(r,r,r)$ for any positive odd integer $r$ in \cite{tsumura2003}. In \cite{tsumura2004}, he obtained the evaluation formula for
$S(p,q,r)$ with $p,q,r\in\mathbb{N}$ and $p+q+r$ odd. To evaluate $R(p,q,r)$, H. Tsumura introduced the partial Tornheim's double series
defined by
\begin{equation}
\mathfrak{T}_{b_1,b_2}(p,q,r):=\sum\limits_{m,n=0}^\infty
\frac{1}{(2m+b_1)^p(2n+b_2)^q(2m+2n+b_1+b_2)^r},
\label{Eq:defn-partial-Tornheim}
\end{equation}
where $b_1,b_2\in\{1,2\}$. Then in \cite[Theorem 4.1]{tsumura2009}, he proved that for any $p,q,r\in\mathbb{N}$ with $r\geqslant 2$
and $p+q+r$ odd, and for
$b_1,b_2\in\{1,2\}$, the values $R(p,q,r)$ and $\mathfrak{T}_{b_1,b_2}(p,q,r)$ can be expressed as polynomials in Riemann zeta values
$\zeta(j)$ ($2\leqslant j\leqslant p+q+r$) with rational coefficients.

In this paper, we give new proofs of two functional relations for $S(s_1,s_2,s_3)$ and $R(s_1,s_2,s_3)$ in Theorem \ref{Thm:funcEq-SR},
from which we obtain new proofs for formulas of $S(p,q,r)$ and $R(p,q,r)$ mentioned in the last paragraph.
The method used here is different from that of H. Tsumura \cite{tsumura2007} and of T. Nakamura \cite{nakamura} for Tornheim's double zeta function case. In fact, it is also valid for proving T. Nakamura's functional relation for Tornheim's double zeta function. For proving this functional relation, we use the following simple facts:
\begin{enumerate}
  \item There is a recursive formula for $T(s_1,s_2,s_3)$:
  $$T(s_1,s_2,s_3)=T(s_1-1,s_2,s_3+1)+T(s_1,s_2-1,s_3+1).$$
  From this formula, finally we only need to treat $T(0,s_2,s_3)$, $T(s_1,0,s_3)$ and $T(s_1,s_2,0)$.
  \item From the definition, we know
  $$T(0,s_2,s_3)=\zeta(s_3,s_2),\; T(s_1,0,s_3)=\zeta(s_3,s_1),\; T(s_1,s_2,0)=\zeta(s_1)\zeta(s_2),$$
  where $\zeta(s_1,s_2)$ is the double zeta function.
  \item We have harmonic shuffle product (stuffle product)
  $$\zeta(s_1)\zeta(s_2)=\zeta(s_1,s_2)+\zeta(s_2,s_1)+\zeta(s_1+s_2).$$
\end{enumerate}
Of course, to avoid singularities we need some modification. Similar facts are also valid for the alternating analogues $R(s_1,s_2,s_3)$ and $S(s_1,s_2,s_3)$.

We give the new proof of the functional relation of Tornheim's double zeta function in Section 2. Then in Section 3, we give our new proof for the functional relations stated in Theorem \ref{Thm:funcEq-SR}.
In Section 4, we give new proofs of H. Tsumura's results mentioned above.

\begin{rem}
After the first version of this paper was submitted to arXiv, Professor Takashi Nakamura kindly informed me (\cite{nakamura10}) that the relation \eqref{Eq:funcEq-SRR} is essentially \cite[Theorme 3.2]{nakamura081}, the relation \eqref{Eq:funcEq-RRS} is proved by \cite[Theorem 3.1]{nakamura082}, Proposition \ref{Prop:R-Odd} and Proposition \ref{Prop:S-Odd} are also \cite[Proposition 3.3]{nakamura082} or \cite[Corollary 1 and Proposition 1]{zhou}, and Proposition \ref{Prop:T} is also \cite[Proposition 5.1]{nakamura082}. I would like to express my gratitude to him for indicating these facts.
\end{rem}


\section{A functional relation for Tornheim's double zeta function}

In \cite[Theorem 4.5]{tsumura2007}, H. Tsumura proved the following functional relation for Tornheim's
double zeta function:
\begin{align}
&T(a,b,s)+(-1)^bT(b,a,s)+(-1)^aT(s,a,b)\label{Eq:funcEq-T-Tsumura}\\
=&\;2\sum\limits_{j=0\atop j\equiv a (2)}^{a}(2^{1-a+j}-1)\zeta(a-j)\sum\limits_{l=0}^{j/2}\frac{(\pi i)^{2l}}{(2l)!}\binom{b-1+j-2l}{j-2l}\zeta(b+j+s-2l)\nonumber\\
&-4\sum\limits_{j=0\atop j\equiv a (2)}^a(2^{1-a+j}-1)\zeta(a-j)\sum\limits_{l=0}^{(j-1)/2}\frac{(\pi i)^{2l}}{(2l+1)!}\sum\limits_{k=0\atop k\equiv b(2)}^{b}\zeta(b-k)\nonumber\\
&\;\;\;\times \binom{k-1+j-2l}{j-2l-1}\zeta(k+j+s-2l),
\nonumber
\end{align}
where (2) means $\hspace{-5pt} \mod 2$, and $a,b\in\mathbb{N}\cup \{0\}$, $b\geqslant 2$, $s\in\mathbb{C}$, except for the singular points of
both sides. In \cite{nakamura}, T. Nakamura gave a \lq\lq simpler" version, which can be restated as the following theorem.

\begin{thm}[{\cite[Theorem 1.2]{nakamura}}]\label{Thm:funcEq-T-Nakamura}
For all $a,b\in\mathbb{N}$ and $s\in\mathbb{C}$ except for the singular points, we have
\begin{equation}
T(a,b,s)+(-1)^bT(b,a,s)+(-1)^aT(s,a,b)=2N(a,b,s)+2N(b,a,s),
\label{Eq:funcEq-T-Nakamura}
\end{equation}
where
$$N(a,b,s):=\sum\limits_{j=0}^{a/2}\binom{a+b-2j-1}{b-1}\zeta(2j)\zeta(a+b+s-2j).$$
\end{thm}

In \cite{matsumoto}, K. Matsumoto, T. Nakamura, H. Ochiai and H. Tsumura showed that the right-hand sides of \eqref{Eq:funcEq-T-Tsumura} and \eqref{Eq:funcEq-T-Nakamura}
are the same. In this section, we first restate their proof with a different method to obtain the key formulas used in the proof.
Then we give a new proof of the functional relation \eqref{Eq:funcEq-T-Nakamura}.

Recall that the Bernoulli polynomials $\{B_n(x)\}$ and the Bernoulli numbers $\{B_n\}$ are defined by
$$\frac{te^{xt}}{e^t-1}=\sum\limits_{n=0}^\infty B_n(x)\frac{t^n}{n!},\qquad
\frac{t}{e^t-1}=\sum\limits_{n=0}^\infty B_n\frac{t^n}{n!},$$
respectively. It is known that $B_n(0)=B_n(1)=B_n$ for any $n\neq 1$ and $B_1(0)=B_1=-B_1(1)=-\frac{1}{2}$.
We recall the formulas
\begin{align}
&\zeta(2n)=-\frac{(2\pi i)^{2n}}{2(2n)!}B_{2n},
\label{Eq:zeta2n}\\
&B_n(1/2)=(2^{1-n}-1)B_n,\label{Eq:bernoulli-1/2}\\
&B_n(x+y)=\sum\limits_{k=0}^n\binom{n}{k}B_k(x)y^{n-k}.
\label{Eq:translation}
\end{align}
From the translation formula \eqref{Eq:translation}, we immediately get the following lemma.

\begin{lem}
For any nonnegative integer $n$, we have
\begin{equation}
\sum\limits_{k=0}^n \binom{2n}{2k}B_{2k}(x)y^{2n-2k}=\frac{1}{2}(B_{2n}(x+y)+B_{2n}(x-y)),
\label{Eq:Bernoulli-Even}
\end{equation}
and
\begin{equation}
\sum\limits_{k=0}^n \binom{2n+1}{2k}B_{2k}(x)y^{2n+1-2k}=\frac{1}{2}(B_{2n+1}(x+y)-B_{2n+1}(x-y)).
\label{Eq:Bernoulli-Odd}
\end{equation}
\end{lem}

Using the above lemma, we get the following key formulas for proving the fact that the right-hand side of \eqref{Eq:funcEq-T-Tsumura} equals
that of \eqref{Eq:funcEq-T-Nakamura}.

\begin{lem}
For any nonnegative integer $n$, we have
\begin{equation}
\sum\limits_{k=0}^n(2^{1-2k}-1)\zeta(2k)\frac{(\pi i)^{2n-2k}}{(2n-2k)!}=\zeta(2n),
\label{Eq:zeta-Even}
\end{equation}
and
\begin{equation}
\sum\limits_{k=0}^n(2^{1-2k}-1)\zeta(2k)\frac{(\pi i)^{2n-2k}}{(2n-2k+1)!}=-\frac{1}{2}\delta_{n,0},
\label{Eq:zeta-Odd}
\end{equation}
where $\delta_{ij}$ is Kronecker's delta symbol.
\end{lem}

\proof We get \eqref{Eq:zeta-Even} from \eqref{Eq:zeta2n}, \eqref{Eq:bernoulli-1/2} and \eqref{Eq:Bernoulli-Even},
and get \eqref{Eq:zeta-Odd} from \eqref{Eq:zeta2n}, \eqref{Eq:bernoulli-1/2} and \eqref{Eq:Bernoulli-Odd}.
\qed

We come to prove the fact that the right-hand sides of \eqref{Eq:funcEq-T-Tsumura} and \eqref{Eq:funcEq-T-Nakamura}
are the same. The right-hand side of \eqref{Eq:funcEq-T-Tsumura} is $2R_1+2R_2$, where
$$R_1=\sum\limits_{j=0\atop j\equiv a (2)}^{a}(2^{1-a+j}-1)\zeta(a-j)\sum\limits_{l=0}^{j/2}\frac{(\pi i)^{2l}}{(2l)!}\binom{b-1+j-2l}{j-2l}\zeta(b+j+s-2l),$$
and
\begin{align*}
R_2=&-2\sum\limits_{j=0\atop j\equiv a (2)}^a(2^{1-a+j}-1)\zeta(a-j)\sum\limits_{l=0}^{(j-1)/2}\frac{(\pi i)^{2l}}{(2l+1)!}\sum\limits_{p=0\atop p\equiv b(2)}^{b}\zeta(b-p)\\
&\times \binom{p-1+j-2l}{j-2l-1}\zeta(p+j+s-2l).
\end{align*}
We show that $R_1=N(a,b,s)$ and $R_2=N(b,a,s)$. For $R_1$, let $a-j=2k$ and $k+l=n$, we get
$$
R_1=\sum\limits_{k=0}^{a/2}(2^{1-2k}-1)\zeta(2k)\sum\limits_{n=k}^{a/2}\frac{(\pi i)^{2n-2k}}{(2n-2k)!}\binom{a+b-2n-1}{a-2n}\zeta(a+b+s-2n).$$
Changing the order of $n$ and $k$, we get
$$R_1=\sum\limits_{n=0}^{a/2}\left(\sum\limits_{k=0}^{n}(2^{1-2k}-1)\zeta(2k)\frac{(\pi i)^{2n-2k}}{(2n-2k)!}\right)\binom{a+b-2n-1}{a-2n}\zeta(a+b+s-2n),$$
which is $N(a,b,s)$ by \eqref{Eq:zeta-Even}.
Similarly for $R_2$, we have
\begin{align*}
R_2=&-2\sum\limits_{k=0}^{a/2}(2^{1-2k}-1)\zeta(2k)\sum\limits_{n=k}^{(a-1)/2}\frac{(\pi i)^{2n-2k}}{(2n-2k+1)!}\sum\limits_{m=0}^{b/2}\zeta(2m)\\
&\times \binom{a+b-2m-2n-1}{a-2n-1}\zeta(a+b+s-2m-2n)\\
=&-2\sum\limits_{n=0}^{(a-1)/2}\left(\sum\limits_{k=0}^n(2^{1-2k}-1)\zeta(2k)\frac{(\pi i)^{2n-2k}}{(2n-2k+1)!}\right)
\sum\limits_{m=0}^{b/2}\zeta(2m)\\
&\times \binom{a+b-2m-2n-1}{a-2n-1}\zeta(a+b+s-2m-2n),
\end{align*}
which is $N(b,a,s)$ by \eqref{Eq:zeta-Odd}. Hence the right-hand sides of \eqref{Eq:funcEq-T-Tsumura} and
\eqref{Eq:funcEq-T-Nakamura} are the same.

In the rest of this section, we give a new proof of the functional relation \eqref{Eq:funcEq-T-Nakamura}. We set
$$F(a,b,s):=T(a,b,s)+(-1)^bT(b,s,a)+(-1)^aT(s,a,b).$$
It is well-known that Tornheim's double zeta function satisfies the following recursive formula
$$T(s_1,s_2,s_3)=T(s_1-1,s_2,s_3+1)+T(s_1,s_2-1,s_3+1).$$
We find that
\begin{align*}
F(a,b,s)=&(T(a-1,b,s+1)+T(a,b-1,s+1))\\
&+(-1)^b(T(b,s+1,a-1)-T(b-1,s+1,a))\\
&+(-1)^a(-T(s+1,a-1,b)+T(s+1,a,b-1)),
\end{align*}
which is just
\begin{equation}
F(a,b,s)=F(a-1,b,s+1)+F(a,b-1,s+1).
\label{Eq:recursive-T}
\end{equation}

As stated in \cite{huard}, we have the following general lemma.

\begin{lem}[\cite{huard}]\label{Lem:rep-general}
Let $X(s_1,s_2,s_3)$ be a function satisfying the recursive relation
$$X(a,b,s)=X(a-1,b,s+1)+X(a,b-1,s+1).$$
Then for any $a,b\in\mathbb{N}$ and $s\in\mathbb{C}$, we have
\begin{align}
X(a,b,s)=&\sum\limits_{j=1}^a\binom{a+b-j-1}{b-1}X(j,0,a+b+s-j)\label{Eq:rep-general}\\
&+\sum\limits_{j=1}^b\binom{a+b-j-1}{a-1}X(0,j,a+b+s-j).
\nonumber
\end{align}
\end{lem}

One can prove this lemma by induction on $a+b$.

Now we want to apply Lemma \ref{Lem:rep-general} to $F(a,b,s)$. Since the singularities of $T(s_1,s_2,s_3)$ lie on the subsets of $\mathbb{C}^3$ defined by one of the equations $s_1+s_3=1-l$, $s_2+s_3=1-l$ ($l\in\mathbb{N}\cup \{0\}$), or $s_1+s_2+s_3=2$, we use an equivalent form of \eqref{Eq:rep-general}
\begin{align}
X(a,b,s)=&\sum\limits_{j=2}^a\binom{a+b-j-1}{b-1}X(j,0,a+b+s-j)\label{Eq:rep-general-2}\\
&+\sum\limits_{j=2}^b\binom{a+b-j-1}{a-1}X(0,j,a+b+s-j)
\nonumber\\
&+\binom{a+b-2}{a-1}X(1,1,a+b+s-2).
\nonumber
\end{align}

We first compute $F(j,0,a+b+s-j)$
and $F(0,j,a+b+s-j)$ for $j\geqslant 2$. It is easy to see $F(j,0,a+b+s-j)=F(0,j,a+b+s-j)$. And we have
\begin{align*}
&F(j,0,a+b+s-j)\\
=&T(j,0,a+b+s-j)+T(0,a+b+s-j,j)+(-1)^jT(a+b+s-j,j,0)\\
=&(1+(-1)^j)\zeta(j)\zeta(a+b+s-j)-\zeta(a+b+s).
\end{align*}
Then we get
\begin{align*}
&\sum\limits_{j=2}^a\binom{a+b-j-1}{b-1}F(j,0,a+b+s-j)\\
=&
2\sum\limits_{j=1}^{a/2}\binom{a+b-2j-1}{b-1}\zeta(2j)\zeta(a+b+s-2j)-\sum\limits_{j=2}^a\binom{a+b-j-1}{b-1}\zeta(a+b+s)\\
=&2\sum\limits_{j=1}^{a/2}\binom{a+b-2j-1}{b-1}\zeta(2j)\zeta(a+b+s-2j)+2\binom{a+b-2}{b}\zeta(0)\zeta(a+b+s),
\end{align*}
and similarly
\begin{align*}
&\sum\limits_{j=2}^b\binom{a+b-j-1}{a-1}F(0,j,a+b+s-j)\\
=&2\sum\limits_{j=1}^{b/2}\binom{a+b-2j-1}{a-1}\zeta(2j)\zeta(a+b+s-2j)+2\binom{a+b-2}{a}\zeta(0)\zeta(a+b+s).
\end{align*}
Here we use the fact that $\zeta(0)=-\frac{1}{2}$. Using a similar result as \cite[Corollary 2]{tornheim}, we get
\begin{align*}
&F(1,1,a+b+s-2)=T(1,1,a+b+s-2)-2T(1,a+b+s-2,1)\\
=&-2\zeta(a+b+s)=4\zeta(0)\zeta(a+b+s).
\end{align*}
Combining the above three equations and \eqref{Eq:rep-general-2},
we finish the proof of Theorem \ref{Thm:funcEq-T-Nakamura}.


\section{Functional relations for $S(s_1,s_2,s_3)$ and $R(s_1,s_2,s_3)$}

As in \cite{zhao}, for complex variables $s_1,s_2$ with $\Re(s_1)>1$ and $\Re(s_2)\geqslant 1$, we define
\begin{align*}
&\zeta(\bar{s}_1,\bar{s}_2):=\sum\limits_{m>n>0}\frac{(-1)^{m+n}}{m^{s_1}n^{s_2}},\\
&\zeta(\bar{s}_1,s_2):=\sum\limits_{m>n>0}\frac{(-1)^{m}}{m^{s_1}n^{s_2}},\\
&\zeta(s_1,\bar{s}_2):=\sum\limits_{m>n>0}\frac{(-1)^{n}}{m^{s_1}n^{s_2}},
\end{align*}
For $\Re(s_2)\geqslant 1$, as in \cite[Proposition 1.1]{AK} we define
\begin{align*}
&\zeta(\bar{1},\bar{s}_2):=\lim\limits_{R\rightarrow\infty}\sum\limits_{R>m>n>0}\frac{(-1)^{m+n}}{mn^{s_2}},\\
&\zeta(\bar{1},s_2):=\lim\limits_{R\rightarrow\infty}\sum\limits_{R>m>n>0}\frac{(-1)^m}{mn^{s_2}}.
\end{align*}
We also define
$$\zeta(\bar{s}):=\sum\limits_{m=1}^\infty \frac{(-1)^m}{m^s}$$
for $\Re(s)>0$. It is easy to see that
\begin{align}
&\zeta(\bar{s})\zeta(t)=\zeta(\bar{s},t)+\zeta(t,\bar{s})+\zeta(\overline{s+t}),\quad (\Re(s)\geqslant 1,\Re(t)>1),
\label{Eq:prod-oneBar}\\
&\zeta(\bar{s})\zeta(\bar{t})=\zeta(\bar{s},\bar{t})+\zeta(\bar{t},\bar{s})+\zeta(s+t), \quad (\Re(s)\geqslant 1,\Re(t)\geqslant 1).
\label{Eq:prod-twoBar}
\end{align}

For $a,b\in\mathbb{N}$ and $s\in\mathbb{C}$, we define
\begin{align*}
&F_1(a,b,s):=S(a,b,s)+(-1)^bR(b,s,a)+(-1)^aR(a,s,b),\\
&F_2(a,b,s):=R(a,b,s)+(-1)^bR(s,b,a)+(-1)^aS(a,s,b).
\end{align*}
Similar to  Tornheim's double zeta function, we have the recursive relations:
\begin{align*}
&S(s_1,s_2,s_3)=S(s_1-1,s_2,s_3+1)+S(s_1,s_2-1,s_3+1),\\
&R(s_1,s_2,s_3)=R(s_1-1,s_2,s_3+1)+R(s_1,s_2-1,s_3+1).
\end{align*}
Then we get the following lemma.

\begin{lem}
We have the following recursive relations:
\begin{align}
&F_1(a,b,s)=F_1(a-1,b,s+1)+F_1(a,b-1,s+1),\label{Eq:recursive-SRR}\\
&F_2(a,b,s)=F_2(a-1,b,s+1)+F_2(a,b-1,s+1).\label{Eq:recursive-RRS}
\end{align}
\end{lem}

Before applying Lemma \ref{Lem:rep-general} to $F_1$ and $F_2$, we make some preparations.
Note that from \cite[Theorem 2.1]{nakamura081}, the singularities of $R(s_1,s_2,s_3)$ lie on the subsets of $\mathbb{C}^3$ defined by one of the equations $s_1+s_3=1-l$ ($l\in\mathbb{N}\cup\{0\}$), and there are no singularities of $S(s_1,s_2,s_3)$.

\begin{lem}
We have
\begin{align}
&F_1(j,0,a+b+s-j)=F_1(0,j,a+b+s-j)\label{Eq:F1-j-0}\\
=&(1+(-1)^j)\zeta(j)\zeta(\overline{a+b+s-j})-\zeta(\overline{a+b+s}), \quad (j\geqslant 2).\nonumber\\
&F_2(j,0,a+b+s-j)=(1+(-1)^j)\zeta(\bar{j})\zeta(\overline{a+b+s-j})-\zeta(a+b+s),\label{Eq:F2-j-0}\\
&F_2(0,j,a+b+s-j)=(1+(-1)^j)\zeta(\bar{j})\zeta(a+b+s-j)-\zeta(\overline{a+b+s}).\label{Eq:F2-0-j}
\end{align}
\end{lem}

\proof We get \eqref{Eq:F1-j-0} and \eqref{Eq:F2-0-j} from \eqref{Eq:prod-oneBar}, and get \eqref{Eq:F2-j-0} from \eqref{Eq:prod-twoBar}.
\qed

Now we have the functional relations.

\begin{thm}[\cite{nakamura081,nakamura082}]\label{Thm:funcEq-SR}
For all $a,b\in\mathbb{N}$ and $s\in\mathbb{C}$ except for the singular points, we have
\begin{equation}
S(a,b,s)+(-1)^bR(b,s,a)+(-1)^aR(a,s,b)=2N_1(a,b,s)+2N_1(b,a,s),
\label{Eq:funcEq-SRR}
\end{equation}
and
\begin{equation}
R(a,b,s)+(-1)^bR(s,b,a)+(-1)^aS(a,s,b)=2N_2(a,b,s)+2N_3(b,a,s),
\label{Eq:funcEq-RRS}
\end{equation}
where
\begin{align*}
&N_1(a,b,s):=\sum\limits_{j=0}^{a/2}\binom{a+b-2j-1}{b-1}(2^{2j+1-a-b-s}-1)\zeta(2j)\zeta(a+b+s-2j),\\
&N_2(a,b,s):=\sum\limits_{j=0}^{a/2}\binom{a+b-2j-1}{b-1}(2^{1-2j}-1)(2^{2j+1-a-b-s}-1)\zeta(2j)\zeta(a+b+s-2j),\\
&N_3(a,b,s):=\sum\limits_{j=0}^{a/2}\binom{a+b-2j-1}{b-1}(2^{1-2j}-1)\zeta(2j)\zeta(a+b+s-2j).
\end{align*}
\end{thm}

\proof We get \eqref{Eq:funcEq-RRS} from \eqref{Eq:recursive-RRS}, \eqref{Eq:F2-j-0}, \eqref{Eq:F2-0-j} and Lemma \ref{Lem:rep-general}. For \eqref{Eq:funcEq-SRR}, we need a formula
\begin{align*}
F_1(1,1,a+b+s-2)=&S(1,1,a+b+s-2)-2R(1,a+b+s-2,1)\\
=&-2\zeta(\overline{a+b+s})=4\zeta(0)\zeta(\overline{a+b+s}),
\end{align*}
which is proved similarly as \cite[Corollary 2]{tornheim}. Using this formula, together with \eqref{Eq:recursive-SRR}, \eqref{Eq:F1-j-0} and \eqref{Eq:rep-general-2}, we get \eqref{Eq:funcEq-SRR}.
\qed

Note that
$$N_1(a,b,s)+N_2(a,b,s)+N_3(a,b,s)=(2^{2-a-b-s}-1)N(a,b,s),$$
where $N(a,b,s)$ is defined in Section 2.


\section{Applications of functional relations}

In \cite[Section 3]{nakamura}, T. Nakamura used the functional relation \eqref{Eq:funcEq-T-Nakamura}
to give new proofs of some formulas for the special values of $T(p,q,r)$ with $p,q,r\in\mathbb{N}$.
For example, we have the evaluation formula of $T(p,q,r)$ when $p+q+r$ is odd as in \cite{huard,nakamura}.

\begin{prop}[\cite{huard,nakamura}]
For  $p,q,r\in\mathbb{N}$ with $p+q+r$ odd, we have
$$T(p,q,r)=(-1)^pN(p,r,q)+(-1)^pN(r,p,q)+(-1)^qN(q,r,p)+(-1)^qN(r,q,p).$$
\end{prop}

In this section, we use the functional relations \eqref{Eq:funcEq-SRR} and \eqref{Eq:funcEq-RRS}
to deduce some formulas for the special values of $S(p,q,r)$, $R(p,q,r)$ and $\mathfrak{T}_{b_1,b_2}(p,q,r)$ with $p,q,r\in\mathbb{N}$
and $b_1,b_2\in\{1,2\}$.

Let $a=b=s=r\in\mathbb{N}$ in \eqref{Eq:funcEq-SRR} and \eqref{Eq:funcEq-RRS}, we get
\begin{align}
&S(r,r,r)+2(-1)^rR(r,r,r)=4N_1(r,r,r),\label{Eq:eva-all-SRR}\\
&(1+(-1)^r)R(r,r,r)+(-1)^rS(r,r,r)=2N_2(r,r,r)+2N_3(r,r,r).\label{Eq:eva-all-RRS}
\end{align}

Let $r=2p$ be even in \eqref{Eq:eva-all-SRR} and \eqref{Eq:eva-all-RRS}, we get a formula which was mentioned in \cite[Eq. (4.2)]{tsumura2002}.

\begin{prop}[{\cite[Eq. (4.2)]{tsumura2002}}]
For any $p\in\mathbb{N}$, we have
\begin{align*}
&S(2p,2p,2p)+2R(2p,2p,2p)\\
=&4\sum\limits_{j=0}^p\binom{4p-2j-1}{2p-1}(2^{2j+1-6p}-1)\zeta(2j)\zeta(6p-2j)\\
=&2\sum\limits_{j=0}^p\binom{4p-2j-1}{2p-1}(2^{2-6p}-2^{2j+1-6p})\zeta(2j)\zeta(6p-2j).
\end{align*}
\end{prop}

The above formulas give some relations for Riemann zeta values. For example, taking $p=1$, we get the relation
$7\zeta(6)=4\zeta(2)\zeta(4)$.

Let $r=2p+1$ be odd in \eqref{Eq:eva-all-RRS}, we get the evaluation formula of $S(2p+1,2p+1,2p+1)$ as in \cite{tsumura2002,matsumoto}.

\begin{prop}[\cite{tsumura2002,matsumoto}]
For any $p\in\mathbb{N}\cup\{0\}$, we have
\begin{align*}
&S(2p+1,2p+1,2p+1)\\
=&2^{-6p}\sum\limits_{j=0}^p\binom{4p+1-2j}{2p}(2^{2j-1}-1)\zeta(2j)\zeta(6p+3-2j).
\end{align*}
\end{prop}

Let $r=2p+1$ be odd in \eqref{Eq:eva-all-SRR}. Using the above formula for $S(2p+1,2p+1,2p+1)$, we get
the evaluation formula of $R(2p+1,2p+1,2p+1)$ as in \cite{tsumura2003,matsumoto}.

\begin{prop}[\cite{tsumura2003,matsumoto}]
For any $p\in\mathbb{N}\cup\{0\}$, we have
\begin{align*}
&R(2p+1,2p+1,2p+1)\\
=&2^{-6p-1}\sum\limits_{j=0}^p\binom{4p+1-2j}{2p}(2^{6p+2}-2^{2j-1}-1)\zeta(2j)\zeta(6p+3-2j).
\end{align*}
\end{prop}

Let $a=p$, $b=q$ and $s=r$ in \eqref{Eq:funcEq-SRR}, we get
\begin{equation}
S(p,q,r)+(-1)^qR(q,r,p)+(-1)^pR(p,r,q)=2N_1(p,q,r)+2N_1(q,p,r).
\label{Eq:lemma-odd-SRR}
\end{equation}
Let $a=p$, $b=r$ and $s=q$ in \eqref{Eq:funcEq-RRS}, we get
$$R(p,r,q)+(-1)^rR(q,r,p)+(-1)^pS(p,q,r)=2N_2(p,r,q)+2N_3(r,p,q),$$
which is
\begin{equation}
(-1)^pR(p,r,q)+(-1)^{p+r}R(q,r,p)+S(p,q,r)=2(-1)^p(N_2(p,r,q)+N_3(r,p,q)).
\label{Eq:lemma-odd-RRS}
\end{equation}
The difference of \eqref{Eq:lemma-odd-SRR} and \eqref{Eq:lemma-odd-RRS} gives
\begin{align*}
&((-1)^q-(-1)^{p+r})R(q,r,p)\\
=& 2N_1(p,q,r)+2N_1(q,p,r)-2(-1)^p(N_2(p,r,q)+N_3(r,p,q)),
\end{align*}
which deduces the evaluation formula of $R(p,q,r)$ when $p+q+r$ is odd as in \cite{tsumura2009}.

\begin{prop}[\cite{tsumura2009}]\label{Prop:R-Odd}
For $p,q,r\in\mathbb{N}$ with $p+q+r$ odd, we have
$$R(p,q,r)=(-1)^pN_1(r,p,q)+(-1)^pN_1(p,r,q)+(-1)^qN_2(r,q,p)+(-1)^qN_3(q,r,p).$$
Explicitly, we have
\begin{align*}
&R(p,q,r)\\
=&(-1)^p\sum\limits_{j=0}^{p/2}\binom{p+r-2j-1}{r-1}(2^{2j+1-p-q-r}-1)\zeta(2j)\zeta(p+q+r-2j)\\
&+(-1)^q\sum\limits_{j=0}^{q/2}\binom{q+r-2j-1}{r-1}(2^{1-2j}-1)\zeta(2j)\zeta(p+q+r-2j)\\
&+(-1)^p\sum\limits_{j=0}^{r/2}\binom{p+r-2j-1}{p-1}(2^{2j+1-p-q-r}-1)\zeta(2j)\zeta(p+q+r-2j)\\
&+(-1)^q\sum\limits_{j=0}^{r/2}\binom{q+r-2j-1}{q-1}(2^{1-2j}-1)(2^{2j+1-p-q-r}-1)\zeta(2j)\zeta(p+q+r-2j).
\end{align*}
\end{prop}

With the help of the above proposition and \eqref{Eq:lemma-odd-SRR}, we get the evaluation formula of $S(p,q,r)$ when $p+q+r$ is odd as in \cite{tsumura2004}.

\begin{prop}[\cite{tsumura2004}]\label{Prop:S-Odd}
For $p,q,r\in\mathbb{N}$ with $p+q+r$ odd, we have
$$S(p,q,r)=(-1)^pN_2(p,r,q)+(-1)^qN_2(q,r,p)+(-1)^pN_3(r,p,q)+(-1)^qN_3(r,q,p).$$
More precisely, we have
\begin{align*}
&S(p,q,r)\\
=&(-1)^p\sum\limits_{j=0}^{p/2}\binom{p+r-2j-1}{r-1}(2^{1-2j}-1)(2^{2j+1-p-q-r}-1)\zeta(2j)\zeta(p+q+r-2j)\\
&+(-1)^q\sum\limits_{j=0}^{q/2}\binom{q+r-2j-1}{r-1}(2^{1-2j}-1)(2^{2j+1-p-q-r}-1)\zeta(2j)\zeta(p+q+r-2j)\\
&+(-1)^p\sum\limits_{j=0}^{r/2}\binom{p+r-2j-1}{p-1}(2^{1-2j}-1)\zeta(2j)\zeta(p+q+r-2j)\\
&+(-1)^q\sum\limits_{j=0}^{r/2}\binom{q+r-2j-1}{q-1}(2^{1-2j}-1)\zeta(2j)\zeta(p+q+r-2j).
\end{align*}
\end{prop}

The evaluation formula for $S(p,q,r)$ with $p+q+r$ odd given by H. Tsumura in \cite{tsumura2004} reads
\begin{align*}
&S(p,q,r)
=(-1)^pN_2(p,r,q)+(-1)^qN_2(q,r,p)\\
&-2(-1)^p\sum\limits_{j=0}^{(r-1)/2}\zeta(\overline{2j})\sum\limits_{\rho=0}^{p/2}\zeta(\overline{2\rho})\sum\limits_{\mu=0}^{(p-2\rho-1)/2}
\binom{p+r-2j-2\rho-2\mu-1}{p-2\rho-2\mu-1}\\
&\times \zeta(p+q+r-2j-2\rho-2\mu)\frac{(\pi i)^{2\mu}}{(2\mu+1)!}\\
&-2(-1)^q\sum\limits_{j=0}^{(r-1)/2}\zeta(\overline{2j})\sum\limits_{\rho=0}^{q/2}\zeta(\overline{2\rho})\sum\limits_{\mu=0}^{(q-2\rho-1)/2}
\binom{q+r-2j-2\rho-2\mu-1}{q-2\rho-2\mu-1}\\
&\times \zeta(p+q+r-2j-2\rho-2\mu)\frac{(\pi i)^{2\mu}}{(2\mu+1)!}.
\end{align*}
The third term of the right-hand side of the above equation equals
\begin{align*}
&-2(-1)^p\sum\limits_{j=0}^{(r-1)/2}\zeta(\overline{2j})\sum\limits_{\rho=0}^{p/2}\zeta(\overline{2\rho})\sum\limits_{n=\rho}^{(p-1)/2}
\binom{p+r-2j-2n-1}{p-2n-1}\\
&\times \zeta(p+q+r-2j-2n)\frac{(\pi i)^{2n-2\rho}}{(2n-2\rho+1)!}.
\end{align*}
Changing the order of $\rho$ and $n$, we see that the above formula equals
\begin{align*}
&-2(-1)^p\sum\limits_{j=0}^{(r-1)/2}\zeta(\overline{2j})\sum\limits_{n=0}^{(p-1)/2}\left(\sum\limits_{\rho=0}^n
\zeta(\overline{2\rho})\frac{(\pi i)^{2n-2\rho}}{(2n-2\rho+1)!}\right)
\\
&\times \binom{p+r-2j-2n-1}{p-2n-1}\zeta(p+q+r-2j-2n),
\end{align*}
and using \eqref{Eq:zeta-Odd}, we find that it becomes
$$(-1)^p\sum\limits_{j=0}^{(r-1)/2}\binom{p+r-2j-1}{p-1}\zeta(\overline{2j})\zeta(p+q+r-2j).$$
Hence the formula of H. Tsumura is nothing but
\begin{align*}
S(p,q,r)
=&(-1)^pN_2(p,r,q)+(-1)^qN_2(q,r,p)\\
&+(-1)^p\sum\limits_{j=0}^{(r-1)/2}\binom{p+r-2j-1}{p-1}\zeta(\overline{2j})\zeta(p+q+r-2j)\\
&+(-1)^q\sum\limits_{j=0}^{(r-1)/2}\binom{q+r-2j-1}{q-1}\zeta(\overline{2j})\zeta(p+q+r-2j).
\end{align*}
Now it is easy to see that the formula of H. Tsumura for $S(p,q,r)$ is the same as that given in Proposition \ref{Prop:S-Odd}.

It is obvious that
\begin{align*}
&\mathfrak{T}_{1,2}(p,q,r)=\mathfrak{T}_{2,1}(q,p,r),\;\;\;\mathfrak{T}_{2,2}(p,q,r)=2^{-p-q-r}T(p,q,r),\\
&R(p,q,r)=-\mathfrak{T}_{1,1}(p,q,r)+\mathfrak{T}_{1,2}(p,q,r)-\mathfrak{T}_{2,1}(p,q,r)+\mathfrak{T}_{2,2}(p,q,r),\\
&S(p,q,r)=\mathfrak{T}_{1,1}(p,q,r)-\mathfrak{T}_{1,2}(p,q,r)-\mathfrak{T}_{2,1}(p,q,r)+\mathfrak{T}_{2,2}(p,q,r).
\end{align*}
Thus we get
\begin{align*}
&\mathfrak{T}_{2,1}(p,q,r)=-\frac{1}{2}(R(p,q,r)+S(p,q,r))+\mathfrak{T}_{2,2}(p,q,r),\\
&\mathfrak{T}_{1,1}(p,q,r)=-\frac{1}{2}(R(p,q,r)+R(q,p,r))+\mathfrak{T}_{2,2}(p,q,r).
\end{align*}
Then we obtain the evaluations of $\mathfrak{T}_{b_1,b_2}(p,q,r)$ when $p+q+r$ is odd as in
\cite{tsumura2003,tsumura2009}.

\begin{prop}[\cite{tsumura2003,tsumura2009}]\label{Prop:T}
For $p,q,r\in\mathbb{N}$ with $p+q+r$ odd, we have
\begin{align*}
\mathfrak{T}_{1,1}(p,q,r)=&-\frac{1}{2}\{(-1)^pN_1(r,p,q)+(-1)^pN_1(p,r,q)+(-1)^pN_2(r,p,q)\\
&+(-1)^pN_3(p,r,q)+(-1)^qN_1(r,q,p)+(-1)^qN_1(q,r,p)\\
&+(-1)^qN_2(r,q,p)+(-1)^qN_3(q,r,p)\}+\mathfrak{T}_{2,2}(p,q,r),
\end{align*}
\begin{align*}
\mathfrak{T}_{1,2}(p,q,r)=&\mathfrak{T}_{2,1}(q,p,r)\\
=&-\frac{1}{2}\{(-1)^pN_2(p,r,q)+(-1)^pN_2(r,p,q)+(-1)^pN_3(p,r,q)\\
&+(-1)^pN_3(r,p,q)+(-1)^qN_1(r,q,p)+(-1)^qN_1(q,r,p)\\
&+(-1)^qN_2(q,r,p)+(-1)^qN_3(r,q,p)\}+\mathfrak{T}_{2,2}(p,q,r),
\end{align*}
and
\begin{align*}
\mathfrak{T}_{2,2}(p,q,r)=&2^{-p-q-r}((-1)^pN(p,r,q)+(-1)^pN(r,p,q)\\
&+(-1)^qN(q,r,p)+(-1)^qN(r,q,p)).
\end{align*}
\end{prop}


\noindent {\bf Acknowledgements.} This work was partially supported by the National Natural Science Foundation of
China (Grant No. 11471245) and Shanghai Natural
Science Foundation (grant no. 14ZR1443500). The author thanks the anonymous referee for his/her helpful comments.


\end{document}